\documentclass[11pt]{amsart}
\usepackage{amsmath}
\newtheorem{thm}{Theorem}[section]
\newtheorem{cor}[thm]{Corollary}
\newtheorem{lem}[thm]{Lemma}
\newtheorem{prop}[thm]{Proposition}
\theoremstyle{definition}

\theoremstyle{remark}
\newtheorem{rem}[thm]{\bf Remark}
\newtheorem{exe}[thm]{Example}
\numberwithin{equation}{section}

\begin{document}
\title{Dynamic of abelian subgroups of GL($n$,
$\mathbb{C}$):\\ a structure's Theorem}

\author{Adlene Ayadi and Habib Marzougui}

\address{Habib Marzougui, Department of Mathematics,  Faculty of Sciences of Bizerte,
Zarzouna. 7021. Tunisia; \ Adlene Ayadi, Department of
Mathematics, Faculty of Sciences of Gafsa, Gafsa, Tunisia}

\email{habib.marzouki@fsb.rnu.tn; adleneth@voila.fr}

\thanks{This work is supported by the research unit: syst\`emes dynamiques et combinatoire: 99UR15-15}



\begin{abstract}

In this paper, we characterize the dynamic of every abelian
subgroups $\mathcal{G}$ of GL($n$, $\mathbb{K}$), $\mathbb{K} =
\mathbb{R}$ or $\mathbb{C}$. We show that there exists a $\mathcal{G}$-invariant, dense
open set $U$ in $\mathbb{K}^{n}$ saturated by minimal orbits
with $\mathbb{K}^{n}- U$ a union of at most $n$ $\mathcal{G}$-invariant
vectorial subspaces of $\mathbb{K}^{n}$ of dimension $n-1$ or
$n-2$ on $\mathbb{K}$. As a consequence, $\mathcal{G}$ has height
at most $n$ and in particular it admits a minimal set in
$\mathbb{K}^{n}-\{0\}$.
\end{abstract}

\maketitle

\section{Introduction} Let $\mathbb{K}=\mathbb{R}$ or $\mathbb{C}$, $GL(n, \ \mathbb{K})$ be the group
of all invertible square matrices of order $n\geq 1$ with entries
in $\mathbb{K}$, and let $\mathcal{G}$ be an abelian subgroup of
GL($n$, $\mathbb{K}$). There is a natural linear action $GL($n$, \
\mathbb{K})\times \mathbb{K}^{n}\longrightarrow \mathbb{K}^{n}:
(A, \ v)\longmapsto Av$. For a vector $v\in \mathbb{K}^{n}$, we
consider the orbit of $\mathcal{G}$ through $v$:

$\mathcal{G}(v) = \{Av, \ A\in \mathcal{G}\} \subset \mathbb{\mathbb{K}}^{n}$. A
subset $E \subset \mathbb{K}^{n}$ is called
\emph{$\mathcal{G}$-invariant} if $A(E)\subset E$ for any $A\in
\mathcal{G}$; that is $E$ is a union of orbits.
\medskip

In \cite{msK04}, Kulikov studied the problem of the
 existence of minimal sets in $\mathbb{R}^{n}-\{0\}$. He constructed an example of discrete
 subgroup of $SL(2, \mathbb{R})$ whose linear action on
$\mathbb{R}^{2}$ is without minimal set in $\mathbb{R}^{2} - \{0 \}$.
F. Dal'bot and A.N. Starkov touched in \cite{fD-anS00} the question of the existence of an
infinitely generated subgroup of $SL(2, \ \mathbb{R})$ with all
orbits dense in $\mathbb{R}^{2}$.

In \cite{aA-hM03}, we studied in the viewpoint closure of orbits the dynamic of a class of
abelian subgroups of $GL(n, \ \mathbb{R})$; those containing an
element \emph{$A\in \mathcal{G}$ which satisfies the condition
($\star$)}: all eigenspaces of $A$ are of dimension $1$ on
$\mathbb{C}$.
\medskip

This work considers the general case: the study of the dynamic of every abelian subgroup
of $GL(n, \ \mathbb{K})$. The purpose here is to develop in this general situation a
setting of a structure's Theorem analogous to a structure's Theorem (in \cite{hM-eS03}) for foliations on closed manifolds.
\medskip

Before stating our main results, we introduce the following notions for groups:
\medskip

A $\mathcal{G}$-invariant subset $E$ of $\mathbb{K}^{n}$ is called
\emph{a minimal set} of $\mathcal{G}$ if every orbit contained in
$E$ is dense in it (this definition it equivalent to say that $E$
is closed in $\mathbb{K}^{n}$, non empty, $\mathcal{G}$-invariant
and has no proper subset with these properties). If $V$ is a
$\mathcal{G}$-invariant open set in $\mathbb{K}^{n}$, a
\emph{minimal set in $V$} is a minimal set of $\mathcal{G}$
restricted to $V$. We say that an orbit $O$ in $V$ is
called \emph{minimal in $V$} if $\overline{O}\cap V$ is
a minimal set in $V$.
\medskip

We call \emph{class of an orbit} $L$ of $\mathcal{G}$ the set
$cl(L)$ of orbits $O$ of $\mathcal{G}$ such that $\overline{O} =
\overline{L}$. If $L$ is an orbit which is minimal in a
$\mathcal{G}$-invariant, open set $V$ then $cl(L) = \overline{L}
\cap V$.
\medskip

An orbit $L$ of $\mathcal{G}$ is said to be at \emph{level} $1$
if, $L$ is minimal in $\mathbb{K}^{n} - \{0\}$. Inductively, we
say that $L$ is at \emph{level} $p$, $p\geq 1$ if every orbit $O
\subset \overline{L} - cl(L)$ is at level $< p$ with at least one
orbit at level $k$ for every $k < p$. The upper bound of levels of
orbits of $\mathcal{G}$ is called the \emph{height} of
$\mathcal{G}$, denoted it by $ht(\mathcal{G})$. For example, if
$ht(\mathcal{G})$ is finite, say $p$ this means that $p$ is the
supremum of $k\in \mathbb{N}$ such that there exist orbits
$\gamma_{1}$, $\gamma_{2}$, ..., $\gamma_{k}$ of $\mathcal{G}$
such that $\overline{\gamma_{1}} \subset
\overline{\gamma_{2}}\subset ... \subset \overline{\gamma_{k}}$
with $\overline{\gamma_{i}} \neq \overline{\gamma_{j}}, \ 1\leq
i\neq j\leq k$.
\medskip

The main result of this paper is the following structure's
Theorem:
\bigskip

{\bf Structure's Theorem.} Let $\mathcal{G}$ be an abelian subgroup
of GL($n$, $\mathbb{K}$) ($\mathbb{K} = \mathbb{R}$ or
$\mathbb{C}$). Then there exists a $\mathcal{G}$-invariant, dense
open set $U$ in $\mathbb{K}^{n}$ with the following properties:

i) Every orbit of $U$ is minimal in $U$.

ii) $\mathbb{K}^{n}-U$ is a union of at most $n$  $\mathcal{G}$-invariants
vectorial subspaces of \;$\mathbb{K}^{n}$ of dimension $n-1$ or
$n-2$ on $\mathbb{K}$.
\medskip

\begin{rem}
If \ $\mathcal{G}$ is an abelian subgroup of \;$GL(n, \mathbb{R})$
containing an element \emph{$A\in \mathcal{G}$ which satisfies
the condition ($\star$)} we know that all orbits of $U$ are homeomorphic
(cf. \cite{aA-hM03}). If \ $\mathcal{G}$ does not contain an element which satisfies the
condition ($\star$) and if $n\geq 3$, this property is not true in general: we give a counterexample
for $n = 3,\ 4, \ 5$ (see Examples 6.1, 6.2  and 6.3) of abelian subgroups $\mathcal{G}$ which
contains two non homeomorphic orbits in any open set in $\mathbb{R}^{3}$(resp. $\mathbb{R}^{4}$;
$\mathbb{C}^{5}$). If $n = 2$, the property remains true: a group $\mathcal{G}$
does not contain an element which satisfies the condition ($\star$) is a subgroup of homotheties. Therefore, if $U = \mathbb{R}^{2}-\{0\}$ and
$u, \ v\in U$, there exists $A\in GL(2,\ \mathbb{R})$ such that $Au = v$ and thus
$A(\mathcal{G}(u)) = \mathcal{G}(v)$.
\end{rem}
\bigskip

\begin{cor} Let $\mathcal{G}$ be an abelian subgroup
of $GL(n, \mathbb{K})$. Then $\mathcal{G}$ has height at most $n$.
\end{cor}
\bigskip

On the other hand, if we remove $\{0 \}$ which is fixed by
$\mathcal{G}$ (i.e. $\mathcal{G}(0) = \{0 \})$, does there exist a
minimal set in $\mathbb{K}^{n} - \{0 \}$? we gave in \cite{aA-hM03} a positive answer for this
question which is here a consequence of Corollary 1.2:
\medskip

\begin{cor}(\cite{aA-hM03})
If \; $\mathcal{G}$ is an abelian subgroup of \;$GL(n,
\mathbb{K})$, it admits a minimal set in \;$\mathbb{K}^{n} -
\{0\}$.
\end{cor}
\medskip

\begin{cor} If $\mathcal{G}$ has a locally dense orbit $O$ in $\mathbb{K}^{n}$ and $C$ a connected component of $U$
meeting $O$ then:

i) $O$ is dense in $C$

ii) Every orbit in $U$ meeting $C$ is dense in it.
\end{cor}
\medskip

\begin{cor} Let \;$\mathcal{G}$ be an abelian subgroup of \;$GL(n,
\mathbb{K})$. If $\mathcal{G}$ has a dense orbit in $\mathbb{K}^{n}$ then every orbit in $U$ is
dense in $\mathbb{K}^{n}$.
\end{cor}
\bigskip

This paper is organized as follows: In section $2$ we give some
notations and technical lemmas. In section $3$, we prove the main theorem for a subgroup of
$S_{n}(\mathbb{K})$. The proof of the structure's Theorem is done in section $4$.
In section $5$, we prove Corollaries 1.2, 1.3, 1.4, 1.5. In section $6$, some examples are given.

\section{Notations and Lemmas}
\bigskip

In this paper, we denote by:

- $\mathbb{K} = \mathbb{R} \ or \ \mathbb{C}$.

- $T_{n}(\mathbb{K})$ the subgroup of $GL(n,\ \mathbb{K})$ of lower triangular matrices,

- $S_{n}(\mathbb{K})$ the subgroup of $T_{n}(\mathbb{K})$ of matrices
$B = (b_{i,j})_{1\leq i, j\leq n}$ with $b_{i,i} = \mu_{B} , \ i= 1, ..., n$.

Every element $B\in S_{n}(\mathbb{K})$ is written in the form
$$B = \left(\begin{array}{cccc}
  \mu_{B} & 0 & .. & 0 \\
  b_{2,1} & \mu_{B} & . & . \\
  . & . & . & 0 \\
  b_{n,1} & .. & b_{n,n-1} & \mu_{B}
\end{array}\right) = \left(\begin{array}{cc}
  B^{(1)} & 0 \\
  L_{B} & \mu_{B}
\end{array}\right),$$ with $B^{(1)}=\left(\begin{array}{cccc}
  \mu_{B} & 0 & .. & 0 \\
  b_{2,1} & \mu_{B} & . & . \\
  . & . & . & 0 \\
  b_{n-1,1} & .. & b_{n-1,n-1} & \mu_{B}
\end{array}\right)\in S_{n-1}(\mathbb{K})$, $\mu_{B}\in \mathbb{K}$ and

$L_{B} = (b_{n,1},....,b_{n,n-1})$.

\bigskip

For this section,  $\mathcal{G}$ is an abelian subgroup of
$S_{n}(\mathbb{K})$.

Denote by:

- $G^{(1)} = \{B^{(1)}, \ B\in \mathcal{G}\}$.

- $\mathcal{C} = (e_{1},...., e_{n})$ the canonical basis of
$\mathbb{K}^{n}$, \ $e_{i} = (0,...,0,1,0,...,0)\in
\mathbb{K}^{n}$ ($1$ is the $i^{th}$ coordinate of $e_{i}$).

- $\mathcal{F} = \{(B - \mu_{B}I_{n})e_{i}\in \mathbb{K}^{n}, \ 1\leq i\leq n-1, \
B\in \mathcal{G}\}$

- $I_{n}$ is the identity matrix of $\mathbb{K}^{n}$.

- $rang(\mathcal{F})$ the rang of $\mathcal{F}$. We have
$rang(\mathcal{F})\leq n-1$.

For every $x = (x_{1},..., x_{n})\in \mathbb{K}^{n}$, we let $x^{(1)} =
(x_{1}, ..., x_{n-1})\in \mathbb{K}^{n-1}$ and $\mathcal{F}^{(1)} = \{(B^{(1)}- \mu_{B}I_{n-1})e^{(1)}_{k},\ \ 1\leq k\leq n-2,\ \ B\in \mathcal{G} \}$.
 We have $x = (x^{(1)}, x_{n})$ and $e_{k} = (e^{(1)}_{k}, 0)$ , $k = 1,...,n-2.$
\bigskip

We start with the following lemmas:
\medskip

\begin{lem} Under the notation above, if $rang(\mathcal{F}) = n-1$ and if
\ $v_{1}, \ ..., \ v_{n-1}$ generate $\mathcal{F}$, then \ $rang(\mathcal{F}^{(1)}) = n-2$ with
$(n-2)$ vectors of $(v^{(1)}_{1},..,v^{(1)}_{n-1})$ generate $\mathcal{F}$.
\end{lem}

\begin{proof}

Let $A_{1}, ..., A_{n-1}\in
\mathcal{G}$ such that $v_{k} = (A_{k}- \mu_{A_{k}}I_{n})e_{i_{k}}$, $k = 1, ..., n-1.$

We let \ $A_{k} = \left(\begin{array}{cc}
  A^{(1)}_{k} & 0 \\
  L_{k} & \mu_{A_{k}}
\end{array}\right)$ and \ $v^{(1)}_{k} = (A^{(1)}_{k}-\mu_{A_{k}}I_{n-1})e^{(1)}_{i_{k}}$.
Then \ $v_{k} = (v^{(1)}_{k}, L_{k}e^{(1)}_{i_{k}})\in vect(e_{2}, ..., e_{n})$ and $v^{(1)}_{k}\in vect(e_{2},...,e_{n-1})$
where $vect(e_{2}, ..., e_{n})$ (resp. $vect(e_{2}, ..., e_{n-1}))$ is the
vectorial subspace generated by $e_{2}, ..., e_{n}$ (resp. $e_{2},...,e_{n-1}$).
Thus, $rang(v^{(1)}_{1}, v^{(1)}_{2}, ..., v^{(1)}_{n-1}) = r \leq n-2$.
If $r < n-2$ there exist $1\leq k \leq n-2$ and
$\alpha_{1},...,\alpha_{k-1},\alpha_{k+1},....\alpha_{n-1}\in\mathbb{K}$ such that
$v_{k}^{(1)} = \alpha_{1}v_{1}^{(1)}
+ ... + \alpha_{k-1}v_{k-1}^{(1)} + \alpha_{k+1}v_{k-1}^{(1)}+
... + \alpha_{n-1}v_{n-1}^{(1)}$.

Lets to prove that $\{
v_{1}^{(1)},..,v_{k-1}^{(1)},v_{k+1}^{(1)},..,v_{n-1}^{(1)}\}$ are linearly independents:
Suppose the contrary ; that is there exists $1\leq
j< k\leq n-1$ such that

$v_{j}^{(1)} = \underset{i\neq k, i\neq j
}{\underset{1\leq i\leq n-1}{\sum}}\beta_{i}v_{i}^{(1)}$. Then :
$v_{j}- \underset{i\neq k, i\neq j }{\underset{1\leq i\leq
n-1}{\sum}}\beta_{i}v_{i} = \beta e_{n}$ and $v_{k}-\underset{i\neq
k}{\underset{1\leq i\leq n-2}{\sum}}\alpha_{i}v_{i} = \alpha e_{n},$
where $\beta = L_{j}e_{i_{j}}^{(1)}- \underset{s\neq k, s\neq j
}{\underset{1\leq s\leq n-2}{\sum}}\beta_{s}L_{s}e_{i_{s}}^{(1)}$
and

$\alpha = L_{k}e_{i_{k}}^{(1)}- \underset{s\neq k}{\underset{1\leq s\leq
n-2}{\sum}}\alpha_{s}L_{s}e_{i_{s}}^{(1)}$. Thus, $$rang(v_{1},...,v_{n-1}) = rang(v_{1},..,v_{j-1},\beta
e_{n},v_{j+1},..,v_{k-1}, \alpha
e_{n},v_{k+1},..,v_{n-1}) < n-1.$$ Therefore $rang(\mathcal{F})= rang(v_{1},...,v_{n-1}) < n-1$, a
contradiction. We deduce that $r = n-2$ and then $rang(\mathcal{F}^{(1)}) = n-2$.
\end{proof}

\bigskip

\begin{lem} Let $u, \ v\in \mathbb{K}^{\star}\times \mathbb{K}^{n-1}$ and
let $B_{m} = \left(\begin{array}{cc}
  B^{(1)}_{m} & 0 \\
  L_{m} & \mu_{m}
\end{array}\right)\in \mathcal{G}$, $m\in \mathbb{N}$,
$B^{(1)}_{m}\in S_{n-1}(\mathbb{K})$ such that
$\underset{m\rightarrow +\infty}{lim}B_{m}u = v$. Suppose that $rang(\mathcal{F}) = n-1$. If
$(B^{(1)}_{m})_{m\in \mathbb{N}}$ is bounded then
$(B_{m})_{m\in \mathbb{N}}$ is bounded.
\end{lem}

\begin{proof}

It suffices to show that $(L_{m})_{m\in \mathbb{N}}$ is bounded.

By Lemma 2.1, it follows that $rang(\mathcal{F}^{(1)}) = n-2$ and there exist $(n-2)$ vectors
of $v_{1}^{(1)}$,..,$v_{n-1}^{(1)}$ which generate $\mathcal{F}$, say
$v_{1}^{(1)}$,...,$v_{n-2}^{(1)}$. We have
$v_{k} = (A_{k}- \lambda_{k}I_{n} )e_{i_{k}}$, \
$v^{(1)}_{k} = (A^{(1)}_{k}- \lambda_{k}I_{n-1} )e^{(1)}_{i_{k}}$,\ \
and $A_{k} = \left(\begin{array}{cc}
  A^{(1)}_{k} & 0 \\
  T_{k} & \lambda_{k}
\end{array}\right)\in \mathcal{G}$, $k = 1, ..., n-2$. Since $A_{k}B_{m} = B_{m}A_{k}$ then
$$L_{m}(A^{(1)}_{k}- \lambda_{k}I_{n-1}) = T_{k}(B^{(1)}_{m}- \mu_{m}I_{n-1})$$

Take $\alpha_{k,m} = T_{k}(B^{(1)}_{m} - \mu_{m}I_{n-1})e_{i_{k}}$,
$k = 1,..., n-2$ and $m\in \mathbb{N}$.

By hypothesis, $(B^{(1)}_{m})_{m\in \mathbb{N}}$ is bounded then
$(B^{(1)}_{m}-\mu_{m}I_{n-1})_{m\in \mathbb{N}}$ is bounded and
therefore $(\alpha_{k,m})_{m\in \mathbb{N}}$ is
bounded for every $k = 1,..., n-2$.

In other part, we have : $$\left\{ \begin{array}{c}
  L_{m}(A^{(1)}_{1}-\lambda_{1}I_{n-1})e_{i_{1}} = \alpha_{1,m} \\
  .................................. \\
  .................................. \\
 L_{m}(A^{(1)}_{n-2}-\lambda_{n-2}I_{n-1})e_{i_{n-2}} = \alpha_{n-2,m}
\end{array}\right.$$

Take $L_{m} = (b_{1,m},..., b_{n-1,m})$ and
$v^{(1)}_{k} = (A^{(1)}_{k}- \lambda_{k}I_{n-1})e_{i_{k}} = (0, a_{k,2},..., a_{k,n-1})$.
Thus, $$\left\{\begin{array}{c}
  a_{1,2}b_{2,m}+ ... +a_{1,n-1}b_{n-1,m} = \alpha_{1,m} \\
  .................................. \\
  .................................. \\
  a_{n-2,2}b_{2,m} + ... + a_{n-2,n-1}b_{n-1,m} = \alpha_{n-2,m}
\end{array}\right.$$
\medskip

This system can be written in the form
$MX_{m}=Y_{m}$, with $X_{m} = \left(\begin{array}{c}
  b_{2,m} \\
  . \\
  . \\
  b_{n-1,m}
\end{array}\right),$
$Y_{m} = \left(\begin{array}{c}
  \alpha_{1,m} \\
  . \\
  . \\
  \alpha_{2,m}
\end{array}\right)$ and $M = \left(\begin{array}{cccc}
  a_{1,2} & .. & .. & a_{1,n-1} \\
  .. & .. &.. & .. \\
  .. & .. & .. & .. \\
  a_{n-2,2} & .. & .. & a_{n-2,n-1}
\end{array}\right)$.

Since $( v^{(1)}_{k})_{1\leq k \leq n-2}$ are independents then
$M$ is invertible, so $M^{-1}Y_{m} = X_{m}$. Since
$(\alpha_{k,m})_{m\in\mathbb{N}}$ is bounded,
$(Y_{m})_{m\in\mathbb{N}}$ is bounded and therefore
$(X_{m})_{m\in\mathbb{N}}$ is bounded i.e.
$(b_{k,m})_{m\in\mathbb{N}}$ is bounded for $k = 2,.., n-1$.

It remains to prove that $(b_{1,m})_{m\in\mathbb{N}}$ is bounded :

Since $\underset{m\rightarrow +\infty}{lim}B_{m}u = v$ we have
$$\underset{m\rightarrow
+\infty}{lim}(b_{1,m}u_{1} + ... + b_{n-1,m}u_{n-1} + \mu_{m}u_{n}) = v_{n}$$
where $u = (u_{1},...,u_{n})$ and $v = (v_{1},...,v_{n})\in
\mathbb{K}^{\star}\times \mathbb{K}^{n-1}$. Or $u_{1}\neq 0$ then $$\underset{m\rightarrow
+\infty}{lim}b_{1,m} = \underset{m\rightarrow
+\infty}{lim}\frac{1}{u_{1}}(v_{n}-b_{2,m}u_{2} - ... -
b_{n-1,m}u_{n-1}-\mu_{m}u_{n}).$$

Since $(b_{k,m})_{m\in \mathbb{N}}$ is bounded for $k=2,..,n-1$,
then $(b_{1,m})_{m\in \mathbb{N}}$ is also bounded. We deduce that
$(L_{m})_{m\in \mathbb{N}}$ is bounded.
\end{proof}

\begin{lem}

Let $r = rang(\mathcal{F})$, $1\leq r\leq n-1$ and
$v_{1},...,v_{r}$ are the generator of $\mathcal{F}$. For every $u\in \mathbb{K}^{n}$,
let $H_{u}$ be the vectorial subspace of $\mathbb{K}^{n}$ generated by $u, v_{1},...,
v_{r}$. Then $H_{u}$ is $\mathcal{G}$-invariant.

In particular, the subspace $F = H_{0}$ generated by $v_{1},...,v_{r}$ is
$\mathcal{G}$-invariant.
\end{lem}

\begin{proof}

Let $w = (z_{1},...,z_{n}) \in H_{u}$ and $B\in \mathcal{G}$ with eigenvalue $\mu$.
We have $Bw = \mu w + (B- \mu I_{n})w$ and $(B- \mu I_{n})w = \underset{i=1}{\overset{n-1}{\sum}}z_{i}(B- \mu
I_{n})e_{i}.$

Since $v_{1},...,v_{r}$ generate $\mathcal{F}$, then for
$k=1,..,n-1$, we have $(B-\mu
I_{n})e_{k}=\underset{i=1}{\overset{r}{\sum}}\beta_{k,i}v_{i}\in
H_{u}.$ Therefore: $Bw = \mu w +\underset{k=1}{\overset{n-1}{\sum}}z_{k}\underset{i=1}{\overset{r}{\sum}}\beta_{k,i}v_{i}\
\ \ \in H_{u}.$
\end{proof}

\section{Proof of the structure's Theorem for subgroups of $S_{n}(\mathbb{K})$}

Let $\mathcal{G}$ be an abelian subgroup of $S_{n}(\mathbb{K})$ and $H$ a
$\mathcal{G}$-invariant vectorial subspace of $\mathbb{K}^{n}$. Recall that
$\mathcal{C} = (e_{1}, ..., e_{n})$ is the canonical basis of $\mathbb{K}^{n}$. Let
 $B\in \mathcal{G}$ and $\varphi_{B}$ the automorphism of $\mathbb{K}^{n}$ with matrix in
$\mathcal{C}$ is $B$. Let $\mathcal{C}_{H}$ be a basis of $H$ and Denote by
$B^{H}$ the matrix of the automorphism restriction $(\varphi_{B})_{/H}$ in $\mathcal{C}_{H}$ and $\mathcal{G}_{H} = \{B^{H},\ \ B\in
\mathcal{G}\}$.

The main result of this section is the following:

\begin{prop} Let $\mathcal{G}$ be an abelian subgroup of
$S_{n}(\mathbb{K})$. For every $u,\
v\in\mathbb{K}^{\star}\times\mathbb{K}^{n-1}$ and for every sequence
$(B_{m})_{m\in \mathbb{N}}$ of $\mathcal{G}$ such that
$\underset{m\rightarrow +\infty}{lim}B_{m}u = v$, there exist a
$\mathcal{G}$-invariant subspace $H$ of $\mathbb{K}^{n}$ and a
basis $\mathcal{C}_{H} = \{u, w_{1},..., w_{p}\}$ of $H$ such
that:

i) $\mathcal{G}_{H}$ is a subgroup of
$S_{p+1}(\mathbb{K})$.

ii) $((B_{m})^{H})_{m\in \mathbb{N}}$ is bounded.
\end{prop}
\bigskip

\begin{rem} In the proposition above, the restriction to a $\mathcal{G}$-invariant
vectorial subspace $H$ is necessary as shown in
Exemple 6.4: there exists a subgroup $\mathcal{G}$ of $GL(4,
\mathbb{R})$ and an umbounded sequence $(\mathcal{B}_{m})_{m\in \mathbb{N}}$
in $\mathcal{G}$ with $\underset{m\rightarrow+\infty}{lim}B_{m}u = v$, for  $u,\ v\in
\mathbb{R}^{\star}\times \mathbb{R}^{3}$.
\end{rem}
\bigskip

\begin{proof}
The proof is by induction on $n$.

For $n = 1$, we have $B_{m} = \lambda_{m}\in \mathbb{K}$ and $u, \ v\in
\mathbb{K}^{\star}$. The condition $\underset{m\rightarrow +\infty}{lim}B_{m}u = v$ shows that
$(\lambda_{m})_{m\in \mathbb{N}}$ is bounded. Then $i)$ and $ii)$ are satisfied for
 $H = \mathbb{K}$ and $\mathcal{E} =(u)$.

Suppose the proposition is true until the order $n-1$.

Let $\mathcal{G}$ be an abelian subgroup of $S_{n}(\mathbb{K})$.
Let $u,\ v\in \mathbb{K}^{\star}\times\mathbb{K}^{n-1}$ and
$(B_{m})_{m\in \mathbb{N}}$ a sequence of $\mathcal{G}$ such that
$\underset{m\rightarrow +\infty}{lim}B_{m}u = v$.

Every $B\in \mathcal{G}$ is written as $$B = \left(\begin{array}{cc}
  B^{(1)} & 0 \\
  L_{B} & \mu_{B}
\end{array}\right),$$ where $B^{(1)}\in S_{n-1}(\mathbb{K})$. Denote by
 $B_{m} = \left(\begin{array}{cc}
  B^{(1)}_{m} & 0 \\
  L_{m} & \mu_{m}
\end{array}\right),$ with $B^{(1)}_{m}\in S_{n-1}(\mathbb{K})$,
$L_{m} = (b^{m}_{n,1}, ..., b^{m}_{n,n-1})$
and by $\mathcal{G}^{(1)} = \{B^{(1)},\ \ B\in \mathcal{G} \}$.
One can check that $\mathcal{G}^{(1)}$ is an abelian subgroup of
$S_{n-1}(\mathbb{K})$.

For every $z = (z_{1}, ..., z_{n})\in \mathbb{K}^{n}$, denote by
$z^{(1)} = (z_{1}, ..., z_{n-1})$.

We let $u = (x_{1}, ..., x_{n})\in
\mathbb{K}^{\star}\times\mathbb{K}^{n-1}$, $v =(y_{1}, ..., y_{n})\in
\mathbb{K}^{n}$, $u^{(1)} = (x_{1}, ..., x_{n-1})$ and $v^{(1)} =
(y_{1}, ..., y_{n-1})$.

As $\underset{m\rightarrow +\infty}{lim}B_{m}u = v$ then
$\underset{m\rightarrow +\infty}{lim}B^{(1)}_{m}u^{(1)} =
v^{(1)}$.

By recurrence hypothesis applied to $\mathcal{G}^{(1)}$ on
$\mathbb{K}^{n-1}$, there exists a $\mathcal{G}^{(1)}$-invariant
vectorial subspace $H_{1}$ of $\mathbb{K}^{n-1}$ and a basis
$\mathcal{C}_{H_{1}} = (u^{(1)},w^{(1)}_{1},...,w^{(1)}_{p})$ of $H_{1}$ such that:

$i)$ $\mathcal{G}^{(1)}_{H_{1}}$ is a subgroup of
$S_{p+1}(\mathbb{K})$ .

ii) $((B^{(1)}_{m})^{H_{1}})_{m\in\mathbb{N}}$ is
bounded.
\bigskip

We distinguish two cases:
\medskip

{\bf {Case 1: $dim(H_{1}) < n-1$.}}

We let $H' = \{ z = (z^{(1)},z_{n}): z^{(1)}\in H_{1},\ \
z_{n}\in \mathbb{K}\}$.

$H'$ is a $\mathcal{G}$-invariant vectorial subspace of $\mathbb{K}^{n}$: if $z = (z^{(1)}, z_{n})\in H'$
and $B\in \mathcal{G}$, then $Bz = \left(\begin{array}{c}
 B^{(1)}z^{(1)}, \ L_{B}z^{(1)}+ \mu_{B}z_{n}
\end{array}\right).$
As, $B^{(1)}z^{(1)}\in H_{1}$ then $Bz\in H'$.

Let $w_{k} = (w^{(1)}_{k}, 0)$, $k = 1,..,p$ and we let $\mathcal{C}_{H'}
= (u, w_{1},.., w_{p}, e_{n})$. $\mathcal{C}_{H'}$ is
a basis of $H'$: if $\alpha, \alpha_{1}, ..., \alpha_{p}, \beta\in \mathbb{K}$ such that
$\alpha u+\beta e_{n}+\underset{i=1}{\overset{p}{\sum}}\alpha_{i}w_{i}=0$ then
$\alpha u^{(1)}+\underset{i=1}{\overset{p}{\sum}}\alpha_{i}w^{(1)}_{i}=0,$
so, $\alpha = \alpha_{1} = ... = \alpha_{p}=0$ and thus
$\beta =0$. Therefore $dim(H') = 1 + dim(H_{1})< n$.
\medskip

Denote by $\mathcal{G}_{H'} = \{B^{H'},\ \ B\in
\mathcal{G}\}$. We will to show that $\mathcal{G}_{H'}$ is a subgroup of $S_{p+2}(\mathbb{K})$:
indeed, if $B\in \mathcal{G}$ and $k = 1,...p$, we have
$B^{(1)}w^{(1)}_{k} = \underset{i=k}{\overset{p}{\sum}}\alpha_{k,i}w^{(1)}_{i}$. So,
$Bw_{k} = (B^{(1)}w^{(1)}_{k}, \ L_{B}w^{(1)}_{k}) = (\underset{i=k}{\overset{p}{\sum}}
\alpha_{k,i}w^{(1)}_{i}, \ L_{B}w^{(1)}_{k})$
$ = (L_{B}w^{(1)}_{k})e_{n}+ \underset{i=k}
{\overset{p}{\sum}}\alpha_{k,i}w_{i}.$ Moreover, $Be_{n} = \mu_{B}e_{n}$. It follows that $B^{H'}\in S_{p+2}(\mathbb{K})$.

In the basis $\mathcal{C}_{H'}$, $u$ (resp. $v$) has coordinate $u_{H'} = (1,0,...,0)$ (resp.
$v_{H'} = (v_{1},...,v_{p+2})$).
We have $\underset{m\rightarrow+\infty}{lim}B_{m}^{H'}u_{H'} = v_{H'}$. So, $\underset{m\rightarrow+\infty}{lim}\mu_{m} =
v_{1}$. As $\underset{m\rightarrow+\infty}{lim}B_{m}u = v$ then
$\underset{m\rightarrow+\infty}{lim}\mu_{m} = \frac{y_{1}}{x_{1}}\neq 0$. It follows that
$v_{1}\neq 0$ and therefore $u_{H'}, \ v_{H'}\in \mathbb{K}^{\star}\times \mathbb{K}^{p+1}$.
\medskip

We can apply the recurrence hypothesis to $\mathcal{G}_{H'}$ on
$\mathbb{K}^{p+2}$, so there exists a
$\mathcal{G}_{H'}$-invariant vectorial subspace $H"$ and a basis
$\mathcal{C}_{H"} = (u_{H'},
w"_{1},..., w"_{q})$ of $H"$ which satisfies the assertions $i)$ and $ii)$.
As $H"$ is $\mathcal{G}_{H'}$-invariant then it is a fortiori
$\mathcal{G}$-invariant. Indeed if $B\in \mathcal{G}$, $v\in H"$ and $\varphi_{B}$ the
automorphism of $\mathbb{K}^{n}$ with matrix in $\mathcal{C}$ is $B$, then $\varphi(v) =
(\varphi_{/ H'})_{/H"}(v)\in H"$ with matrix in $\mathcal{C}$ is $(B^{H'})^{H"}$.
\bigskip

{\bf {Case 2: $dim(H_{1}) = n-1$}}

In this case, $H_{1} = \mathbb{K}^{n-1}$ and then $(B^{(1)}_{m})^{H_{1}} =
B^{(1)}_{m}$ is bounded by hypothesis.
\medskip

One distinguish three cases:
\medskip

Let $\mathcal{F} = \{(B- \mu_{B}I_{n})e_{k},\ 1\leq k\leq n-1, \ B\in \mathcal{G}\}$ and
$r = rang(\mathcal{F})$.
\medskip

{\bf {Case 2. a):  $r = n-1$}}

By Lemma 2.2, $(B_{m})_{m\in \mathbb{N}}$ is bounded. The assertions i) and ii) of the
proposition follow by taking $H = \mathbb{K}^{n}$ and $\mathcal{C} = (u, e_{2}, ..., e_{n})$
a basis of $\mathbb{K}^{n}$.
\medskip

{\bf {Case 2. b): $1\leq r < n-1$ }}

Let $v_{1}, ..., v_{r}$ generate $\mathcal{F}$. Let $H_{u}$(resp. $F$) be the vectorial subspace of $\mathbb{K}^{n}$
generated by $(u, v_{1}, ..., v_{r})$ (resp. $(v_{1}, ..., v_{r})$). By Lemma 2.3, $H_{u}$
(resp. $F$) are
$\mathcal{G}$-invariant. Let $(w^{(2)}_{1}, ..., w^{(2)}_{r})$ be a basis of $F$ such that
for every $B\in \mathcal{G}$, $B^{F}$ is lower triangular. Let $w_{k} = (0,w^{(2)}_{k}),
\ k = 1,..., r$. Then $\mathcal{C}_{H_{u}} = (u,
w_{1}, ..., w_{r})$ is a basis of $H_{u}$. For every $B\in \mathcal{G}$, $B^{H_{u}}$ is
lower triangular. Then $\mathcal{G}_{H_{u}}$ is a subgroup of $S_{r+1}(\mathbb{K})$.
Moreover, $u, \ v\in \mathbb{K}^{\star}\times \mathbb{K}^{r}$.
Since $dim(H_{u}) = r+1 < n$, then the proposition follows by applying the recurrence
hypothesis on $\mathcal{G}_{H_{u}}$.
\bigskip

{\bf {Case 2. c): $r = 0$}}

In this case, for every $B\in \mathcal{G}$, $(B- \mu_{B}I_{n})e_{k} = 0$, $k =1, ..., n$. Then,
$B = \mu_{B}I_{n}$ and therefore $\mathcal{G}$ is an abelian
subgroup of homotheties of $GL(n, \mathbb{K})$. By taking $H = \mathbb{K}^{n}$ and
$\mathcal{C} = (u, e_{2},..., e_{n})$ a basis of $\mathbb{K}^{n}$, the assertion i) and ii) of proposition are satisfied:

We have $B_{m} = \mu_{m}I_{n}$, $m\in \mathbb{N}$ and
$\underset{m\rightarrow +\infty}{lim}\lambda_{m}u = v$ then $\underset{m\rightarrow
+\infty}{lim}\mu_{m} = \frac{y_{1}}{x_{1}}$. Hence $(B_{m})_{m\in \mathbb{N}}$ is bounded.
\end{proof}
\bigskip

\begin{cor} Let $\mathcal{G}$ be an abelian subgroup of $S_{n}(\mathbb{K})$
($\mathbb{K} = \mathbb{R}$ or $\mathbb{C}$). For every $u,\ v\in
\mathbb{K}^{\star}\times \mathbb{K}^{n-1}$ and for every sequence
$(B_{m})_{m\in \mathbb{N}}$ of $\mathcal{G}$ such that
$\underset{m\rightarrow +\infty}{lim}B_{m}u = v$, we have
$\underset{m\rightarrow +\infty}{lim}B_{m}^{-1}v = u.$
\end{cor}
\medskip

\begin{proof}

Let $u,\ v\in \mathbb{K}^{n}$ and $(B_{m})_{m\in \mathbb{N}}$ a
sequence in $\mathcal{G}$ such that $\underset{m\rightarrow
+\infty}{lim}B_{m}u = v$. By Proposition 3.1, there exists a
$\mathcal{G}$-invariant vectorial subspace $H$ of $\mathbb{K}^{n}$
containing $u,\ v$ and a basis $\mathcal{C}_{H} = (u, w_{1}, ..., w_{p})$ of $H$
such that $(B_{m}^{H})_{m\in \mathbb{N}}$ is bounded.

Denote by $u_{H} = (u_{1}, ..., u_{n})$ (resp. $v_{H} = (v_{1}, ..., v_{n})$) where $u_{1}, ..., u_{n}$
(resp. $v_{1}, ..., v_{n}$) are the coordinate of $u$ (resp. $v$) in $\mathcal{C}_{H}$. Denote by $\parallel \ \parallel$
the norm on $GL(n, \mathbb{K})$ defined by : $\parallel B(x) = \underset{x\in \mathbb{K}^{n} -\{0\}}{sup}
(\frac{\parallel Bx \parallel}{\parallel x \parallel})$.
We have
$$\parallel(B^{H}_{m})^{-1}v_{H} - u_{H}\parallel \leq \parallel
(B^{H}_{m})^{-1}\parallel \ \parallel v_{H} - B^{H}_{m}u_{H}
\parallel.$$

The matrix $N_{m} = B^{H}_{m} - \mu_{m}I_{p+1}$ is nilpotent of
order $p+1$. Then

$(B^{H}_{m})^{-1} = \frac{1}{\mu_{m}}\underset{k=0}{\overset{p}{\sum}}
(-1)^{k}\frac{1}{(\mu_{m})^{k}}N_{m}^{k}$. Therefore :

$$\parallel (B^{H}_{m})^{-1}\parallel
\leq \underset{k=0}{\overset{p}{\sum}}\frac{1}{|\mu_{m}|^{k+1}}\parallel
N_{m}\|^{k}$$

Since $(B_{m}^{H})_{m\in \mathbb{N}}$ is bounded and
$\underset{m\rightarrow+\infty}{lim}B_{m}u = v$, $u,\
v\in \mathbb{K}^{\star}\times \mathbb{K}^{n-1}$, then
$(\frac{1}{\mu_{m}})_{m\in\mathbb{N}}$ is bounded. Therefore, $(N_{m})_{m\in
\mathbb{N}}$ and then $((B^{H}_{m})^{-1})_{m\in\mathbb{N}}$ is bounded. Since $\underset{m\rightarrow +\infty}{lim}B^{H}_{m}u_{H} =
v_{H}$, we deduce that
$\underset{m\rightarrow+\infty}{lim}(B_{m}^{H})^{-1}v_{H} = u_{H}$. Hence,
$\underset{m\rightarrow +\infty}{lim}B^{-1}_{m}v = u$.
\end{proof}

\section{Proof of the structure's Theorem : General Case}

In fact, we will prove a slightly strong result; that is:
\bigskip

\begin{thm} Let $\mathcal{G}$ be an abelian subgroup of $GL(n, \mathbb{K})$ and let $E_{i} = \{x = (x_{1}, ..., x_{n})
\in \mathbb{K}^{n}: x_{i} = 0\}$. Then there exists $P\in GL(n, \
\mathbb{K})$ and finitely many $\mathcal{G}$-invariant subspaces $H_{k}$ with
$H_{k} = P(E_{i_{k}})$ if $\mathbb{K} = \mathbb{C}$, and $P(E_{i_{k}})$ or
$P(E_{i_{k}}\cap E_{i_{k}+1})$,
if $\mathbb{K} = \mathbb{R}$, $1\leq k \leq r$, $1\leq i_{1}, ...,
i_{r}\leq n$ such that $\mathcal{G}$ satisfies the property $\mathcal{P}$:

if $u,\ v\in U = \mathbb{K}^{n}
-\underset{k=1}{\overset{r}{\cup}}H_{k}$ and $(B_{m})_{m\in
\mathbb{N}}$ is a sequence of $\mathcal{G}$ such that
$\underset{m\rightarrow +\infty }{lim}B_{m}u = v$ then
$\underset{m\rightarrow +\infty}{lim}B^{-1}_{m}v = u$.
\end{thm}
\bigskip

The proof of the structure's Theorem is completed as follows:

From Theorem 4.1, we let $U = \mathbb{K}^{n}
-\underset{k=1}{\overset{r}{\cup}}H_{k}$. It is clear that $U$ is
a $\mathcal{G}$-invariant dense open set in $\mathbb{K}^{n}$ and
the property $\mathcal{P}$ implies in particular that every orbit
of $U$ is minimal in $U$.
\bigskip

\begin{center}{\bf Case $\mathbb{K} = \mathbb{C}$} \end{center}
\medskip

The proof uses induction on $n$.

The Theorem is true for $n=1$: take $P = I_{\mathbb{C}}$, $H_{1} = \{0\}$ and
$U = \mathbb{C}^{\star}$. So the property $\mathcal{P}$ is satisfied.

Suppose the theorem is true until the order $n-1$. Let $\mathcal{G}$ be an abelian subgroup of $GL(n, \mathbb{C})$.
We distinguish two cases:
\medskip

{\emph{Case 1: Every element of $\mathcal{G}$ has only one eigenvalue}}
\medskip

In this case, there exists $P\in GL(n, \mathbb{C})$ such that
$\mathcal{G}' = P^{-1}\mathcal{G}P$ is a subgroup of
$S_{n}(\mathbb{C})$. By taking $H_{1} =
P(E_{1})$ and $U = \mathbb{C}^{n} - H_{1} =
P\left(\mathbb{C}^{\star}\times\mathbb{C}^{n-1}\right)$, then, by Corollary $3.3$, the property
$\mathcal{P}$ follows.
\medskip

{\emph{Case 2: There exists $A\in \mathcal{G}$ having at least two complex eigenvalues}}

In this case if $\lambda_{1}$, ..., $\lambda_{p}$ be the
eigenvalues of $A$ with order of multiplicities $n_{1}, ...,
n_{p}$ respectively and $E_{k} = Ker(A- \lambda_{k}I_{n})^{n_{k}}$
be the characteristic space of $A$ associated to $\lambda_{k}$ then $1\leq n_{k} < n$.
The space $E_{k}$ ($1\leq k \leq p$) is $\mathcal{G}$-invariant: indeed, if
$B\in G$ and $x\in E_{k}$ then $(A- \lambda_{k}I_{n})^{n_{k}}B(x) =
B(A- \lambda_{k}I_{n})^{n_{k}}(x) = 0$. Denote by $\mathcal{G}_{E_{k}} = \{B^{E_{k}},\
B\in \mathcal{G}\}$. Then $\mathcal{G}_{E_{k}}$ is an abelian
subgroup of $GL(n_{k}, \mathbb{C})$, $k =1,..,p$.

Using the recurrence hypothesis on $\mathcal{G}_{E_{k}}$, there exist
$P_{k}\in GL(n_{k}, \mathbb{C})$ and finitely many subspaces
$H^{k}_{1} = P_{k}(E_{k,j_{1}}), ..., H^{k}_{r_{k}} = P_{k}(E_{k,j_{r_{k}}})\subset E_{k}$, where
$1\leq j_{1}, ..., j_{r_{k}}\leq n_{k}$, and for $i = 1, ..., r_{k}$, $E_{k, j_{i}} = \{x = (x_{k, 1}, ..., x_{k, n_{k}})
\in \mathbb{C}^{n_{k}}: x_{k, j_{i}} =0\}$ such that the
property $\mathcal{P}$ is satisfied; if $u_{k}, \ v_{k}\in U_{k} =
\mathbb{C}^{n_{k}} - \underset{i=1}{\overset{r_{k}}{\cup}}H^{k}_{i}$ and
$(B^{(k)}_{m})_{m\in \mathbb{N}}$ be a sequence in
$\mathcal{G}_{k}$ such that
$\underset{m\rightarrow+\infty}{lim}B^{(k)}_{m}u_{k}=v_{k}$ then
$\underset{m\rightarrow+\infty}{lim}(B^{(k)}_{m})^{-1}v_{k}=u_{k}$.

Since $dim(E_{k}) = n_{k}$, $k = 1,..,p$, and
$\underset{k=1}{\overset{p}{\bigoplus}}E_{k} = \mathbb{C}^{n}$, we let
$\mathcal{C}_{k} = (e_{1,k}, ..., e_{1,n_{k}}$) a basis of
$E_{k}$, $k = 1,.., p$. Hence,
$\mathcal{C'} = \underset{k=1}{\overset{p}{\cup}}\mathcal{C}_{k}$
is a basis of $\mathbb{C}^{n}$. Let $Q\in GL(n, \mathbb{C})$ the
matrix of base change from $\mathcal{C}$ to $\mathcal{C'}$. Then
for every $B\in \mathcal{G}$, we have:

$$Q^{-1}BQ = \left(\begin{array}{ccc}
  B^{E_{1}} & 0 & 0 \\
  0 & . & 0 \\
  0 & 0 & B^{E_{p}}
\end{array}\right)$$

where $B^{E_{k}}\in GL(n_{k}, \mathbb{C})$. Take $P =
\left(\begin{array}{ccc}
  P_{1} & 0 & 0 \\
  0 & . & 0 \\
  0 & 0 & P_{p}
\end{array}\right)$, $i = 1,..,r_{k}$, $k = 1,..,p$. We let $F_{k, j_{i}} = \{x =
(x_{1,1}, ...,
x_{1,n_{1}}; ...; x_{p, 1}, ..., x_{p, n_{p}}) \in \mathbb{C}^{n}:
x_{k, j_{i}} = 0\}$, $K_{i}^{k} = QP(F_{k,j_{i}})$ and
$U = \mathbb{C}^{n}- \underset{k=1}{\overset{p}{\bigcup}}(\underset{i=1}{\overset{r_{k}}
{\bigcup}}K^{k}_{i})$.

Let $u, \ v \in U$ and a sequence $(B_{m})_{m\in \mathbb{N}}\subset
\mathcal{G}$ such that $\underset{m\rightarrow+\infty}{lim}B_{m}u =
v$. Then $\underset{m\rightarrow +\infty}{lim} Q^{-1}B_{m}Q(Q^{-1}u) =
Q^{-1}v$. Take $B'_{m} = Q^{-1}B_{m}Q = \left(\begin{array}{ccc}
  B^{E_{1}}_{m} & 0 & 0 \\
  0 & . & 0 \\
  0 & 0 & B^{E_{p}}_{m}
\end{array}\right)$, $u' = Q^{-1}u = u'_{1}+ ...+ u'_{p}$ and $v' = Q^{-1}v =
v'_{1}+....+ v'_{p}$, where $u'_{k}, v'_{k}\in E_{k},\ \ k = 1,...,p$. We have
$u', \ v'\in U = \mathbb{C}^{n}- \underset{k=1}{\overset{p}{\bigcup}}\left(\underset{i=1}
{\overset{r_{k}}{\bigcup}}P(F_{k,j_{i}})\right)$. Hence, $u'_{k}, v'_{k}\in U_{k}$,
$k =1,...,p$. From $\underset{m\rightarrow +\infty}{lim}B_{m}u = v$ we have
$\underset{m\rightarrow +\infty}{lim}B'_{m}u' = v'$. So,
$\underset{m\rightarrow +\infty}{lim}B^{E_{k}}_{m}u'_{k} = v'_{k}$, $k =1,...,p$.
Therefore $\underset{m\rightarrow
+\infty}{lim}\left(B^{E_{k}}_{m}\right)^{-1}v'_{k} = u'_{k}$, for
every $k = 1,...,p$.

Since $(B'_{m})^{-1} = \left(\begin{array}{ccc}
  \left(B^{E_{1}}_{m}\right)^{-1} & 0 & 0 \\
  0 & . & 0 \\
  0 & 0 & \left(B^{E_{p}}_{m}\right)^{-1}
\end{array}\right)$, we deduce that

$\underset{m\rightarrow +\infty}{lim}\left(B'_{m}\right)^{-1}v' = u'$ and therefore
$\underset{m\rightarrow +\infty}{lim}\left(B_{m}\right)^{-1}v = u$.
\medskip

\begin{center}{\bf Case $\mathbb{K} = \mathbb{R}$} \end{center}

The proof uses the following lemma :
\begin{lem} if $\mathcal{G}$ is an abelian subgroup of $GL(n, \mathbb{R})$ which contains
a matrix $A$ having only two conjugates complex eigenvalues $\lambda$ and
$\overline{\lambda}$, then there exists $P\in GL(n, \
\mathbb{R})$ and finitely many $\mathcal{G}$-invariant subspaces
$H_{k} = P(E_{i_{k}}\cap E_{i_{k}+1})$, $1\leq k \leq r$, $1\leq i_{1}, ...,
i_{r}\leq n$, such that $\mathcal{G}$ satisfies the property $\mathcal{P}$.
\end{lem}
\medskip

\begin{proof}
Denote by $E_{\lambda} = Ker(A- \lambda I_{n})^{s}\subset
\mathbb{C}^{n}$(resp. $E_{\overline{\lambda}} =
Ker(A- \overline{\lambda} I_{n})^{s}\subset \mathbb{C}^{n}$) the
characteristic space of $A$ associated to $\lambda$ (resp.
$\overline{\lambda}$). Then $E_{\lambda}$ and $E_{\overline{\lambda}}$
are $\mathcal{G}$-invariant and $E_{\overline{\lambda}} =
\overline{E_{\lambda}}$.

Take $(v_{1}, ..., v_{s})$ be a basis of $E_{\lambda}$. Hence,
$\mathcal{E} = (v_{1},..., v_{s}, \overline{v_{1}},...,
\overline{v_{s}})$ is a basis of $\mathbb{C}^{n}$. Every matrix $B\in\mathcal{G}$ is written in $\mathcal{E}$ in the
form :
$$B = \left(\begin{array}{cc}
    B^{E_{\lambda}} & 0 \\
    0 & \overline{B^{E_{\lambda}}} \\
    \end{array}\right).$$

Denote by $\mathcal{G}_{E_{\lambda}} = \{B^{E_{\lambda}},\ \
B\in\mathcal{G}\}$. $\mathcal{G}_{E_{\lambda}}$ is an abelian
subgroup of $GL(s, \mathbb{C})$. By applying the case $\mathbb{K}
= \mathbb{C}$, there exists $P_{s}\in GL(s, \ \mathbb{C})$ and
finitely many $\mathcal{G}_{E_{\lambda}}$-invariant hyperplans
$H^{E_{\lambda}}_{k} = P_{s}(E_{i_{k}})_{1\leq k \leq r}$, $1\leq
i_{1}, ..., i_{r}\leq s$, where $E_{i} = \{z = (z_{1}, ...,z_{s})
\in \mathbb{C}^{s}: z_{i} = 0\}$, such that $\mathcal{G}_{E_{\lambda}}$ satisfies the property
$\mathcal{P}$ on $U_{s} = \mathbb{C}^{s} - \underset{k=1}{\overset{r}{\bigcup}}H^{E_{\lambda}}
_{k}.$ Denote by  $P = \left(\begin{array}{cc}
  P_{s} & 0 \\
  0 & \overline{P_{s}}
\end{array}\right)\in GL(n, \mathbb{C})$ and
$\mathcal{E}' = P^{-1}(\mathcal{E}) = (w_{1},..,w_{s},\overline{w_{1}},..,\overline{w_{s}})$,
 where  $w_{k} = P^{-1}(v_{k}) = P_{s}^{-1}(v_{k})$ and  $\overline{w_{k}} =
 P^{-1}(\overline{v_{k}}) = \overline{P_{s}^{-1}}(\overline{v_{k}})$,
 $k = 1,..,s$.
\medskip

Take $e'_{2k-1} = \frac{w_{k} +
\overline{w_{k}}}{2}$, $e'_{2k} = \frac{w_{k} -
\overline{w_{k}}}{2i}$, $k = 1,..,s$. Then
$\mathcal{C}' = (e'_{1},..., e'_{n})$ is a basis of
$\mathbb{R}^{n}$. Denote by $Q\in GL(n, \mathbb{R})$ be the matrix of
basis change of $\mathcal{C}$ to $\mathcal{C}'$.
For $i = 1,..., s$, we let $F_{i} = \{x = (x_{1},y_{1}, ...,
x_{s},y_{s}) \in \mathbb{R}^{n}: x_{i} = y_{i}= 0\}$. Hence,
$F_{i} = E_{2i-1}\cap E_{2i}$. Take $H_{k} = Q(F_{i_{k}})$ and $U =
\mathbb{R}^{n}- \underset{k=1}{\overset{r}{\cup}}H_{k}$, $k = 1,...,r$.
Lets show the property $\mathcal{P}$ for $\mathcal{G}$ on $U$:

Let $u = (u_{1},...,u_{n}), \ v =(v_{1},...,v_{n})\in U$ and
$(B_{m})_{m\in \mathbb{N}}\subset \mathcal{G}$ such that
$\underset{m\rightarrow +\infty }{lim}B_{m}u = v$.
Denote by $Q^{-1}u=(x_{1},y_{1},.....,x_{s},y_{s}),\ Q^{-1}v=(x'_{1},y'_{1},...,x'_{s},y'_{s})$. Since $\mathcal{C}'$ is a basis of
$\mathbb{C}^{n}$, we let $R\in GL(n, \mathbb{C})$ the matrix of
change of $\mathcal{C}'$ to $\mathcal{E}'$. Hence $u' = R^{-1}Q^{-1}u
= (z_{1}, ..., z_{s},\overline{z_{1}},..., \overline{z_{s}})$ and
$v' = R^{-1}Q^{-1}v = (z'_{1}, ..., z'_{s},\overline{z'_{1}}, ...,
\overline{z'_{s}})$, where $z_{k} = \frac{x_{k}- iy_{k}}{2}$,
$z'_{k} = \frac{x'_{k}- iy'_{k}}{2}$.

Since $Q^{-1}u, \ Q^{-1}v\in
Q^{-1}(U) = \mathbb{R}^{n}- \underset{k=1}{\overset{r}{\bigcup}}F_{i_{k}}$
then $z_{i_{k}}, \ z'_{i_{k}} \in \mathbb{R}^{\star}$. It follows that
$(z_{1},...,z_{s})\in\mathbb{C}^{s}- \underset{k=1}{\overset{r}{\bigcup}}E_{i_{k}} =
P_{s}^{-1}(U_{s})$, and
$(\overline{z_{1}},...,\overline{z_{s}})\in \overline{P_{s}}^{-1}(\overline{U_{s}})$. Hence,
$u', \ v' \in P_{s}^{-1}(U_{s})\times \overline{P_{s}}^{-1}(\overline{U_{s}}) =
P^{-1}(U_{s}\times \overline{U_{s}})$. Then,
 $$ (1) \ Pu', \ Pv' \in U_{s}\times \overline{U_{s}}.$$

From $\underset{m\rightarrow +\infty}{lim}B_{m}u = v$, we have

$$(2) \ \underset{m\rightarrow +\infty}{lim}(S^{-1}B_{m}SPu' = Pv'.$$

where $S = QRP^{-1}$. As $S$ is the matrix of basis change of $\mathcal{C}$
to $\mathcal{E}$ then, $S^{-1}B_{m}S$ is the matrix of $B_{m}$
in the basis $\mathcal{E}$. Hence :

  $$S^{-1}B_{m}S = \left(\begin{array}{cc}
    (B_{m})^{E_{\lambda}} & 0 \\
    0 & \overline{(B_{m})^{E_{\lambda}}} \\
  \end{array}\right)$$

  So, $(1)$ and $(2)$ imply that $\underset{m\rightarrow
  +\infty}{lim}B_{m}^{E_{\lambda}}(z_{1},..,z_{s})=(z'_{1},..,z'_{s}).$

By applying the property $\mathcal{P}$ on $U_{s}$, we obtain :

$\underset{m\rightarrow +\infty}{lim}(B_{m}^{E_{\lambda}})^{-1}(z'_{1},..,z'_{s}) =
(z_{1},..,z_{s})$ and $\underset{m\rightarrow +\infty}{lim}(\overline{B_{m}^{E_{\lambda}}})^{-1}
(\overline{z'_{1}},..,\overline{z'_{s}}) = (\overline{z_{1}},..,\overline{z_{s}}).$

Hence, $\underset{m\rightarrow +\infty}{lim}(S^{-1}B^{-1}_{m}S)Pv' = Pu'$ and therefore
$\underset{m\rightarrow +\infty}{lim}B^{-1}_{m}v = u $.
\end{proof}
\medskip

\begin{proof} of Theorem 4.1. Case: $\mathbb{K} = \mathbb{R}$

The proof is by induction on $m$.

The case $n = 1$ is the same as for $\mathbb{K} = \mathbb{C}$.

Suppose the Theorem is true until the order $n-1$. We distinguish two cases:
\bigskip

{\emph{Case 1: Every element of $\mathcal{G}$ has only one real
eigenvalue or only two non real conjugates complex eigenvalues}}

One distinguish two cases:
\medskip

{\emph{Case 1. a): there exist $A\in \mathcal{G}$ which has
only two conjugates complex eigenvalues $\lambda$ and $\overline{\lambda}$}}.

The result follows in this case from Lemma 4.2.
\medskip

{\emph {Case 1. b): Every matrix of $\mathcal{G}$ has only one real
eigenvalue}}.

There exists $Q\in GL(n, \mathbb{R})$ such that
$Q\mathcal{G}Q^{-1}$ is a subgroup of $S_{n}(\mathbb{R})$. We
conclude as for $\mathbb{K} = \mathbb{C}$ in the case $1$.
\medskip

{\emph {Case 2: there exists $A\in \mathcal{G}$ having at least two
non conjugate complex eigenspaces $\lambda$ and $\mu$ such that
$\lambda\neq\overline{\mu}$}}.
\medskip

Let $\lambda_{1},..., \lambda_{r}$ be the complex eigenvalues
of $A$ in $\mathbb{C}- \mathbb{R}$ and its conjugate of
multiplicities $m_{1},..., m_{r}$ and $\alpha_{1},...,
\alpha_{p}$ the real eigenvalues of $A$ with multiplicities
$n_{1},..., n_{p}$ respectively.

Denote by

- $F_{l} = Ker((A- \lambda_{l}I_{n})(A- \overline{\lambda_{l}}I_{n}))^{m_{l}}$,
$l =1,...,r$.

- $E_{k} = Ker(A- \alpha_{k}I_{n})^{n_{k}}$ the characteristic space of
$A$ associated to $\alpha_{k}$, $k = 1,...,p$.

We have $\underset{k=1}{\overset{p}{\bigoplus}}E_{k}\oplus
\underset{l=1}{\overset{r}{\bigoplus}}F_{l} = \mathbb{R}^{n}$,
$dim(E_{k}) = n_{k}$ and $dim(F_{l}) = 2m_{l}$,
$k = 1,..,p$, $l = 1,..,r$. The case $2$ implies that $n_{k}< n$ and $2m_{l}< n$.

The spaces $E_{k}$ and $F_{l}$ are $\mathcal{G}$-invariant: if $B\in \mathcal{G}$ and
$x\in F_{l}$ then :

$((A- \lambda_{l}I_{n})(A- \overline{\lambda_{l}}I_{n}))^{m_{l}}B(x) = B((A- \lambda_{l}I_{n})
(A- \overline{\lambda_{l}}I_{n}))^{m_{l}}(x) = 0$, $l= 1,...,r$
and $(A- \alpha_{k}I_{n})^{m_{k}}B(x) = B(A- \alpha_{k}I_{n})^{n_{k}}(x) = 0$, $k = 1,...,p$.

Let $\mathcal{G}_{E_{k}} = \{ M^{E_{k}},\ M\in \mathcal{G}\}$ and
$\mathcal{G}_{F_{l}} = \{ M^{F_{l}},\ M\in \mathcal{G}\}$.
Then $\mathcal{G}_{E_{k}}$ and $\mathcal{G}_{F_{l}}$ are abelian
subgroup of $GL(n_{k}, \mathbb{R})$ and $GL(2m_{l}, \mathbb{R})$ respectively.

We get the result by applying for $\mathcal{G}_{E_{k}}$ and $\mathcal{G}_{F_{l}}$ the same
proof as for $\mathbb{K} = \mathbb{C}$ in the case $2$.
\end{proof}

\section{Proof of corollaries }

{\bf Proof of Corollary $1.2$.} We process inductively on $n$. If
$n = 1$, $\mathcal{G}$ is a subgroup of homotheties of
\;$\mathbb{K}$. Then it is clear that every orbit $O$ in
$\mathbb{K} -\{0\}$ is minimal in it. Therefore, \;$O$\; is at
level \;$1$. Now, let $L$ be an orbit of $\mathcal{G}$ in
$\mathbb{K}^{n}$. If $L$ is contained in one of the $n$ subspaces, say
$H$ and since $H$ is $\mathcal{G}$-invariant and of dimension
$n-1$ or $n-2$ then by inductive hypothesis, $L$ is at level
$p \leq n-1$ or $n-2$.

If $L$ is contained in $U$, then by the structure's Theorem, $L$
is minimal in $U$. It follows that $cl(L) = \overline{L} \cap U$.
So, for every orbit $O\subset \overline{L} - cl(L)$ we have
$O\subset \mathbb{R}^{n} - U$. Then $O$ is in one of the subspaces
of dimension $n-1$ or $n-2$. By inductive hypothesis applied
to the restriction $\mathcal{G}_{H}$ of $\mathcal{G}$
 to $H$, $O$ is at level
$\leq n-1$ or $n-2$. it follows that $L$ is at level $\leq n$.
$\Box$. \;
\bigskip

{\bf Proof\ of\ Corollary\ $1.3$.}
Take an orbit closure in $\mathbb{K}^{n} - \{0\}$. If it happens not to be minimal in
$\mathbb{K}^{n} - \{0\}$, take a smaller one. If it happens not to be minimal, take a smaller one, we can repeat
 similar until the $n$th step, then by Corollary 1.2, the $n$th one is necessarily minimal.
\bigskip

{\bf Proof\ of\ Corollary\ $1.4$.}

Let $O$ be a dense orbit in $\mathbb{K}^{n}$ (i.e. $\overline{O} = \mathbb{K}^{n}$). Then
$O\subset U$ and $\overline{O}\cap U = U$. Since $O$ is minimal in $U$ then for every orbit
$L\subset U$, we have $\overline{L}\cap U = \overline{O}\cap U$. Therefore
$\overline{L} = \overline{O} = \mathbb{K}^{n}$. \ \ \ \ \ \ \ \ \
\ \ \ \ \ \ \ \ \ \ \ \ \ \  \ \  \ \ \  \ \ \ \ \ \ \ \ \ \ \ \ \
\ \ \ \ \ \ \ \ \ \ \  \ \  \ \ \  \ \ \ \ \ \ \ \ \ \ \ \ \ \ \ \
\ \ \ \ \ \ \ \  \ \  \ \ \ \ $\Box$
\medskip

{\bf Proof\ of\ Corollary\ $1.5$.} \;

If $O$ is a locally dense orbit in $\mathbb{K}^{n}$ (
i.e. $\overset{\circ}{\overline{O}}\neq \emptyset$) then
$O\subset U$. Let $\mathcal{C}$ be a connected component
of $U$ meeting $O$. Then $\overline{O}\cap
\mathcal{C}$ is a non empty closed subset in $\mathcal{C}$. Lets
show that $\overline{O}\cap \mathcal{C}$ is open in
$\mathcal{C}$. Let $v\in\overline{O}\cap \mathcal{C}$. Since
$O$ is minimal in $U$ then $\overline{O}\cap
U = \overline{\mathcal{G}(v)}\cap U$. So,
$\overset{\circ}{\overline{\mathcal{G}(v)}} =
\overset{\circ}{\overline{O}}\neq \emptyset$. Then
$v\in \overset{\circ}{\overline{\mathcal{G}(v)}}\cap
\mathcal{C} \subset\overline{O}\cap \mathcal{C}$.

Lets show that every orbit meeting $\mathcal{C}$ is dense in
$\mathcal{C}$: if $O'$ is an orbit meeting $\mathcal{C}$
then $O'\subset U$. Since $O$ is dense in
$\mathcal{C}$ then $O'\subset \overline{O}$
and then
$\overline{O'}\cap\mathcal{C} = \overline{O}\cap\mathcal{C} = \mathcal{C}$.\
\ \ \ \ \ \ \ \ \ \ \ \ \ \ \ \ \ \ \ \ \ \  \ \  \ \ \  \ \ \ \ \
\ \ \ \ \ \ \ \ \ \ \ \ \ \ \ \ \ \ \  \ \  \ \ \  \ \ \ \ \ \ \ \
\ \ \ \ \ \ \ \ \ \ \ \ \ \ \ \  \ \  \ \ \ \ $\Box$
\bigskip

\section{{\bf  Examples}}

\begin{exe}

Let $\mathcal{G}$ be the subgroup generated by

$A = \left(\begin{array}{ccc}
  1 & 0 & 0 \\
  0 & 1 & 0 \\
  1 & 0 & 1
\end{array}\right)$ and $B = \left(\begin{array}{ccc}
  1 & 0 & 0 \\
  0 & 1 & 0 \\
  0 & 1 & 1
\end{array}\right)$. Then :

i) if $u\in \mathbb{Q}^{\star}\times
\mathbb{Q}\times\mathbb{R}$, $\mathcal{G}(u)$ is closed in $\mathbb{R}^{3}$.

ii) if $u\in \mathbb{Q}^{\star}\times (\mathbb{R}-\mathbb{Q})\times\mathbb{R}$,
$\mathcal{G}(u)$ is dense in a straight line.
\end{exe}

\begin{proof} Let $u = (x,y,z)\in \mathbb{R}^{3}$. We have
$$\mathcal{G}(u) = \{(x, y, nx + my + z): \ \ n,\ m \in \mathbb{Z}^{2}\}.$$

i) if $u\in \mathbb{Q}^{\star}\times \mathbb{Q}\times\mathbb{R}$ then $x\mathbb{Z}+ y\mathbb{Z}$
is closed in $\mathbb{R}$ and therefore $\mathcal{G}(u)$ is closed in $\mathbb{R}^{3}$.

ii) if $u\in \mathbb{Q}^{\star}\times (\mathbb{R}-\mathbb{Q})\times\mathbb{R}$, then
$x\mathbb{Z}+ y\mathbb{Z}$ is dense in $\mathbb{R}$ and therefore $\mathcal{G}(u)$ is
dense in a straight line.
\end{proof}
\medskip

\begin{exe} Let $\mathcal{G}$ be the subgroup of $GL(4, \mathbb{R})$
generated by

$A = \left(\begin{array}{cccc}
  1 & 0 & 0 & 0 \\
  0 & 1 & 0 & 0 \\
  0 & 0 & 1 & 0 \\
  1 & 0 & 0 & 1
\end{array}\right)$ and  $B = \left(\begin{array}{cccc}
  1 & 0 & 0 & 0 \\
  0 & 1 & 0 & 0 \\
  0 & 0 & 1 & 0 \\
  0 & 1 & 0 & 1
\end{array}\right)$

Then:

i) if $u\in (\mathbb{Q}^{\star})^{2}\times \mathbb{R}^{2}$, $\mathcal{G}(u)$ is closed in $\mathbb{R}^{4}$.

ii) if $u\in \mathbb{Q}^{\star}\times(\mathbb{R}-\mathbb{Q})\times \mathbb{R}^{2}$,
$\mathcal{G}(u)$ is dense in a straight line.
\end{exe}

\begin{proof} Take $u = (x, y, z, t)\in \mathbb{R}^{4}$ and $\mathcal{G}(u)$ be the orbit of $u$.
We have $$\mathcal{G}(u) = \{A^{n}B^{m}u: \ \ n, m \in \mathbb{Z}\ \} = \{(x, y, z,
 nx + my+ t): \ \ n, m\in \mathbb{Z}\ \}.$$

i) if $x,\ y\in \mathbb{Q}^{\star}$ then $x\mathbb{Z}+ y\mathbb{Z}$ is closed in
$\mathbb{R}$. So, $\mathcal{G}(u)$ is closed in $\mathbb{R}^{4}$.

ii) if $x\in \mathbb{Q}^{\star}$ and $y\in
\mathbb{R}- \mathbb{Q}$ then  $x\mathbb{Z}+ y\mathbb{Z}$ is dense in
$\mathbb{R}$. So, $\mathcal{G}(u)$ is dense in a straight line of $\mathbb{R}^{4}$.
\end{proof}
\medskip

\begin{exe} Let $\mathcal{G}$ be the abelian subgroup of $GL(5, \mathbb{C})$
generated by

$A = \left(\begin{array}{ccccc}
  1 & 0 & 0 & 0 & 0\\
  0 & 1 & 0 & 0 & 0\\
  0 & 0 & 1 & 0 & 0 \\
  0 & 0 & 0 & 1 & 0 \\
  1 & 0 & 0 & 0 & 1
\end{array}\right)$, $B = \left(\begin{array}{ccccc}
  1 & 0 & 0 & 0 & 0\\
  0 & 1 & 0 & 0 & 0\\
  0 & 0 & 1 & 0 & 0 \\
  0 & 0 & 0 & 1 & 0 \\
  0 & 1 & 0 & 0 & 1
\end{array}\right)$ and $C = \left(\begin{array}{ccccc}
  1 & 0 & 0 & 0 & 0\\
  0 & 1 & 0 & 0 & 0\\
  0 & 0 & 1 & 0 & 0 \\
  0 & 0 & 0 & 1 & 0 \\
  0 & 0 & 1 & 0 & 1
\end{array}\right)$

Then:

 i) if $z\in (\mathbb{Q}^{\star}+i\mathbb{Q}^{\star})^{3}\times \mathbb{C}^{2}$, $\mathcal{G}(z)$ is closed in $\mathbb{C}^{5}$.

ii) if $z = (1+i,\sqrt{3}+i\sqrt{2}, \sqrt{2}+i, z_{4}, z_{5})$ with  $z_{4}, \
z_{5}\in \mathbb{C}^{2}$ then $\mathcal{G}(z)$ is dense in a complex straight line.
\end{exe}

\begin{proof}

Take $z = (z_{1},z_{2},z_{3},z_{4}, z_{5})\in \mathbb{C}^{5}$. We have:
 $$\mathcal{G}(z) = \{A^{n}B^{m}C^{p}z: \ \ n, m, p\in
\mathbb{Z}\} =$$ $$\left \{(z_{1}, z_{2}, z_{3}, z_{4}, nz_{1}+
 mz_{2}+ pz_{3}+ z_{5} \right): \ \ n, m, p\in \mathbb{Z}\}$$

i) Suppose that $z\in (\mathbb{Q}^{\star}+ i \mathbb{Q}^{\star})^{3}\times
\mathbb{C}^{2}$ and $z_{k} = x_{k}+ i y_{k}, \ k = 1, 2, 3$, so
$x_{k}, \ y_{k}\in \mathbb{Q}^{\star}$. The subgroup $H = \mathbb{Z}(x_{1},y_{1})+
\mathbb{Z}(x_{2}, y_{2})+\mathbb{Z}(x_{3},y_{3})$ of $\mathbb{R}^{2}$
is then closed in $\mathbb{R}^{2}$. It follows that \
$z_{1}\mathbb{Z}+ z_{2}\mathbb{Z}+ z_{3}\mathbb{Z}$ \ is closed in $\mathbb{C}$ and
therefore $\mathcal{G}(z)$ is closed in $\mathbb{C}^{5}$.

ii) Lets show first that the subgroup
$H = \mathbb{Z}(1,1)+ \mathbb{Z}(\sqrt{3},\sqrt{2})+\mathbb{Z}(\sqrt{2},1)$ of $\mathbb{R}^{2}$
is dense in $\mathbb{R}^{2}$. By Kronecker generalized Theorem (\cite{mW-cou}), it suffices to show
that if $(s_{1}, s_{2}, s_{3})\in \mathbb{Z}^{3}-\{0\}$,

$$ \Delta = \left|\begin{array}{ccc}
  1 & \sqrt{3} & \sqrt{2} \\
  1 & \sqrt{2} & 1 \\
  s_{1} & s_{2} & s_{3}
\end{array}\right| = \sqrt{3}(s_{1}-s_{3})+\sqrt{2}(s_{3}+s_{2})-(s_{2}+2s_{1})\neq 0:$$
suppose the contrary; that is $\Delta=0$. Since $\sqrt{3}$, $\sqrt{2}$ and $1$ are rationally independents
then $\left\{\begin{array}{c}
  s_{1}-s_{3}=0 \\
  s_{3}+s_{2}=0 \\
  s_{2}+2s_{1}=0
\end{array}\right.$ thus, $s_{1} = s_{2} = s_{3} = 0$, absurd.

if $z = (1+i,\sqrt{3}+i\sqrt{2}, \sqrt{2}+i, z_{4}, z_{5})$ with  $(z_{4}, \
z_{5})\in \mathbb{C}^{2}$ then $\mathcal{G}(z)$ is homeomorphic to $H$, so $\mathcal{G}(z)$ is
dense in a complex straight line of $\mathbb{C}^{5}$.
\end{proof}
\medskip

\begin{exe} Let $\mathcal{G}$ be the subgroup of $GL(4, \mathbb{R})$ generated by

$A = \left(\begin{array}{cccc}
  1 & 0 & 0 & 0 \\
  0 & 1 & 0 & 0 \\
  0 & 0 & 1 & 0 \\
  \sqrt{2}-1 & 1 & 0 & 1
\end{array}\right)$ and $B = \left(\begin{array}{cccc}
  1 & 0 & 0 & 0 \\
  0 & 1 & 0 & 0 \\
  0 & 0 & 1 & 0 \\
  1 & 0 & 0 & 1
\end{array}\right)$, and let

$u = (1,1,0,0); \ v = (1,1,0,\sqrt{3})$. Then:

i) $v\in \overline{\mathcal{G}(u)}- \mathcal{G}(u)$ .

ii) there exists an unbounded sequence $(B_{m})_{m\in \mathbb{N}}$ in $\mathcal{G}$ such that

$\underset{m\rightarrow +\infty}{lim}B_{m}u = v$

iii) The vectorial subspace $H_{u}$ generated by $u$ and $e_{4}=(0,0,0,1)$
is $\mathcal{G}$-invariant and the restriction $(B_{m}^{H_{u}})_{m\in \mathbb{N}}$ is bounded.
\end{exe}

\begin{proof}

i) We have $\mathcal{G}(u) = \{A^{n}B^{m}u: \ \ n, m\in\mathbb{Z}\ \} = \{(1,\ 1,\
0,\ n\sqrt{2} + m ): \ \ n, m\in \mathbb{Z}\ \}$. Since $\mathbb{Z}\sqrt{2}+ \mathbb{Z}$
is dense in $\mathbb{R}$, there exist $k_{m}, \ s_{m}\in \mathbb{Z}$ such that
$\underset{m\rightarrow +\infty}{lim}k_{m}\sqrt{2}+s_{m} = \sqrt{3}$. So,
$\underset{m\rightarrow +\infty}{lim}A^{k_{m}}B^{s_{m}}u = v$ and then
$v\in \overline{\mathcal{G}(u)}$. Since $\sqrt{2}, \ \sqrt{3}$ and $1$ are rationally independents
then $\sqrt{3}\notin \mathbb{Z}\sqrt{2}+ \mathbb{Z}$ and therefore
$v\notin \mathcal{G}(u)$. It follows that we can choose \ $k_{m}, \ s_{m}$ \ such that $\underset{m\rightarrow+\infty}{lim}|k_{m}| = +\infty$.

ii) Take $B_{m} = A^{k_{m}}B^{s_{m}}$. We have
$B_{m} = \left(\begin{array}{cccc}
  1 & 0 & 0 & 0 \\
  0 & 1 & 0 & 0 \\
  0 & 0 & 1 & 0 \\
  k_{m}(\sqrt{2}-1)+s_{m} & k_{m} & 0 & 1
\end{array}\right)$. Since
$\underset{m\rightarrow+\infty}{lim}|k_{m}| = +\infty$ then $(B_{m})_{m\in \mathbb{N}}$ is unbounded.

iii) $H_{u}$ is $\mathcal{G}$-invariant: if $w = \alpha u + \beta e_{4}\in H_{u}$, we have

$Aw = (\alpha, \alpha, 0 ,  \alpha\sqrt{2} + \beta) = \alpha u +
(\alpha\sqrt{2}+\beta)e_{4}\in H_{u}$ and

$Bw =(\alpha, \alpha, 0, \alpha +\beta ) = \alpha u + (\alpha + \beta)e_{4}\in H_{u}$.

Let $A^{H_{u}}$, $B^{H_{u}}$ and $B^{H_{u}}_{m}$ be respectively the restriction of $A$, $B$ and
$B_{m}$ to $H_{u}$ in the basis $(u, e_{4})$. Then
$A^{H_{u}} = \left(\begin{array}{cc}
  1 & 0 \\
  \sqrt{2} & 1
\end{array}\right)$, $B^{H_{u}} = \left(\begin{array}{cc}
  1 & 0 \\
  1 & 1
\end{array}\right)$ and $B^{H_{u}}_{m} = \left(\begin{array}{cc}
  1 & 0 \\
  k_{m}\sqrt{2}+s_{m} & 1
\end{array}\right)$. Since
$\underset{m\rightarrow+\infty}{lim}k_{m}\sqrt{2} + s_{m} = \sqrt{3}$,
 $(B^{H_{u}}_{m})_{m\in \mathbb{N}}$ is bounded.
\end{proof}

\bibliographystyle{amsplain}
\bibliography{xbib}

\vskip 0,4 cm

\end{document}